\documentclass{amsart}
\usepackage{amsaddr}
\usepackage{hyperref}
\usepackage{etoolbox}

\usepackage{graphicx}
\usepackage{tikz}
\usepackage{tikz-cd}
\usepackage{verbatim}
\usepackage{setspace}
\usepackage{mathrsfs}

\newtheorem{theorem}{Theorem}[section]
\newtheorem{lemma}[theorem]{Lemma}

\theoremstyle{definition}
\newtheorem{definition}[theorem]{Definition}

\newtheorem{examples}[theorem]{Examples}

\theoremstyle{remark}

\theoremstyle{proposition}
\newtheorem{proposition}[theorem]{Proposition}

\theoremstyle{corollary}

\numberwithin{equation}{section}

 \DeclareMathOperator{\Ext}{Ext}

\DeclareMathOperator{\E}{E}
\newcommand{\abs}[1]{\lvert#1\rvert}

\newcommand*\xbar[1]{%
\hbox{%
\vbox{%
\hrule height 0.5pt 
\kern0.5ex
\hbox{%
\kern-0.1em
\ensuremath{#1}%
\kern-0.1em
}%
}%
}%
}
\newcommand\restr[2]{{
\left.\kern-\nulldelimiterspace 
#1 
\right|_{#2} 
}}

\begin{document}
\title{On (Co)Pure Baer Injective Modules}

\author{Mohanad Farhan Hamid}
\address{Department of Production and Metallurgy Engineering,
\\University of Technology, Baghdad, Iraq
\\E-mails: {70261@uotechnology.edu.iq
\newline mohanadfhamid@yahoo.com}}

 \curraddr{}

\thanks{}





\maketitle

\section*{\bf{Abstract}}
\noindent For a given class of $R$-modules $\mathcal{Q}$, a module $M$ is
called $\mathcal{Q}$-\emph{copure Baer injective} if any map from a $\mathcal{Q}$-copure
left ideal of $R$ into $M$ can be extended to a map from $R$ into
$M$. Depending on the class $\mathcal{Q}$, this concept is both a dualization and
a generalization of pure Baer injectivity. We show that every
module can be embedded as $\mathcal{Q}$-copure submodule of a $\mathcal{Q}$-copure
Baer injective module. Certain types of rings are characterized using
properties of $\mathcal{Q}$-copure Baer injective modules. For example a ring $R$ is $\mathcal{Q}$-coregular if and only if every $\mathcal{Q}$-copure Baer injective $R$-module is injective.

\bigskip
\noindent \textbf{\keywordsname.} $\mathcal{Q}$-copure submodule, $\mathcal{Q}$-copure
Baer injective module, pure Baer injective module.

\section{\bf{Introduction}}
\noindent Let $\mathcal{Q}$ be a non-empty class of left $R$-modules. An exact sequence
\[
0 \rightarrow A \stackrel{f}\rightarrow B \stackrel{g}\rightarrow C \rightarrow 0 \tag{1} \label{seq}
\]

of left $R$-modules is called $\mathcal{Q}$-\emph{copure} if every module in $\mathcal{Q}$ is injective with respect to the sequence. In this case, $f$ is called a $\mathcal{Q}$-copure monomorphism and $g$ a $\mathcal{Q}$-copure epimorphism \cite[p.322]{Wis}.  If we denote by $\mathcal{PI}$ the class of pure injective modules then the $\mathcal{PI}$-copure sequences are exactly the pure exact ones, see \cite[p.290]{Wis}. So not only does this concept dualize purity  but generalizes it as well. We will need the following lemma later.

\begin{lemma} \label{pushout} \cite[p.323]{Wis} For a given class of modules $\mathcal{Q}$, the following hold.
\begin{enumerate}
\item Any pushout of a $\mathcal{Q}$-copure monomorhphism is a $\mathcal{Q}$-copure monomorphism.
\item If $g\circ f$ in sequence \ref{seq} above is a $\mathcal{Q}$-copure monomorhphism then so is $f$.
\end{enumerate}
\end{lemma}
For details about $\mathcal{Q}$-copure submodules the reader is referred to section 38 of \cite{Wis}.

Thani \cite{Nada} introduced pure Baer injective modules as those modules which are injective with respect to all pure exact sequences with the ring $R$ as a middle term. Here we study $\mathcal{Q}$-copure Baer injective modules for some given non-empty class of left $R$-modules $\mathcal{Q}$, i.e. modules injective with respect to all $\mathcal{Q}$-copure sequences with $R$ as a middle term. Pure Baer injective modules are, now, a special case of $\mathcal{Q}$-copure Baer injectives by choosing $\mathcal{Q} = \mathcal{PI}$.

Unless otherwise stated the ring $R$ is always associative with identity, all modules are left unital $R$-modules, and $\mathcal{Q}$ is a non-empty class of modules. If there is no confusion or if the class $\mathcal{Q}$ is known we will drop the letter $\mathcal{Q}$ and just say copure sequences and copure Baer injective modules.

\section{\bf{Copure Baer Injective Modules}}
\begin{definition}
An $R$-module $M$ is called $\mathcal{Q}$-\emph{copure Baer injective} if any homomorphism from a $\mathcal{Q}$-copure left ideal of $R$ into $M$ has an extension to a homomorphism from $R$ into $M$.
\end{definition}

We will often write copure Baer injective and mean $\mathcal{Q}$-copure Baer injective for some given class $\mathcal{Q}$, just like when we say module and homomorphism (or map) and mean $R$-module and $R$-homomrphism (or $R$-map) for some given ring $R$.

\begin{examples}
\begin{enumerate}
\item Injective modules are $\mathcal{Q}$-copure Baer injective for any class $\mathcal{Q}$.
\item All pure Baer injective (and therefore all pure injective) modules are $\mathcal{PI}$-copure Baer injective.
\item Putting the class $\mathcal{Q}= \{\mathbb{Z}\}$, we see that none of the proper ideals of $\mathbb{Z}$ is $\{\mathbb{Z}\}$-copure. Hence all $\mathbb{Z}$-modules are $\{\mathbb{Z}\}$-copure Baer injective but of course not all of them are injective.

\item We know that all $\mathbb{Z}$-modules are pure Baer injective, however, not all of them are $\mathcal{Q}$-copure Baer injective for all classes $\mathcal{Q}$. For example, let the class $\mathcal{Q} = \{\mathbb{Z}_2\}$. The sequence $0 \rightarrow \mathbb{Z}_3 \rightarrow \mathbb{Z}_9$ is $\mathcal{Q}$-copure exact because the only map $\mathbb{Z}_3 \rightarrow \mathbb{Z}_2$ is the zero map which can, of course, be extended to a map $\mathbb{Z}_9 \rightarrow \mathbb{Z}_2$. But $\mathbb{Z}_3$ is not injective with respect to the above sequence, hence it is not $\mathcal{Q}$-copure Baer injective.

\item Let the ring $R$ be $\mathbb{Z}_4$ and $\mathcal{Q} = \{\mathbb{Z}_4\}$. Since $\mathbb{Z}_4$ is quasi injective, the sequence $0 \rightarrow  \mathbb{Z}_2 \rightarrow \mathbb{Z}_4$ is $\mathcal{Q}$-copure exact. It is in fact the only nontrivial one! So, both of $\mathbb{Z}_4$ and $\mathbb{Z}_3$ are $\mathcal{Q}$-copure Baer injective, while $\mathbb{Z}_2$ is not. To see this consider the following diagram:
$$\begin{tikzcd}
\mathbb{Z}_2 \arrow[hook]{r}\arrow{d}{1_{\mathbb{Z}_2}}
&\mathbb{Z}_4 &\\
\mathbb{Z}_2
\end{tikzcd}$$
which cannot be completed because $\mathbb{Z}_2$ is not a direct summand of $\mathbb{Z}_4$.

\item If $\mathcal{Q}$ is the class of simple modules then the class of copure Baer injective modules equals the class $\mathscr{M}$ of modules injective with respect to all inclusions $I \rightarrow R$ with $I$ an $s$-pure left ideal of $R$, see \cite{ICrivei}.

\item Any module $Q$ is, of course, $\{Q\}$-copure Baer injective but may not, in general, be pure Baer injective. 
\end{enumerate}

\end{examples}

The following proposition is easy to verify.

\begin{proposition} \label{Ext}
\begin{enumerate}
\item The direct product (resp., direct sum) of a (finite) family of modules is copure Baer injective if and only if each factor is copure Baer injective.
\item An $R$-module $M$ is copure Baer injective if and only if $\Ext(R/I,M)=0$ for every copure left ideal $I$ of $R$.
\end{enumerate}
\end{proposition}

\begin{proposition}
The class of copure Baer injective modules is closed under extensions.
\begin{proof}
Let $0 \rightarrow A \rightarrow B \rightarrow C \rightarrow 0$ be an exact sequence with $A$ and $C$ copure Baer injective. Exactness of the sequence $0 \rightarrow \Ext(R/I,A) \rightarrow \Ext(R/I,B) \rightarrow \Ext(R/I,C) \rightarrow 0$ gives, by Proposition \ref{Ext}, that $\Ext(R/I,B) = 0$ for any copure left ideal $I$ of $R$. 
\end{proof}
\end{proposition}

Thani \cite{Nada} introduced left pure hereditary rings as those rings whose every pure left ideal is projective. Here, we define left copure hereditary rings.

\begin{definition} The ring $R$ is called \emph{left $\mathcal{Q}$-copure hereditary} if every  copure left ideal of $R$ is projective.
\end{definition}
 Of course, left pure hereditary rings are $\mathcal{PI}$-copure hereditary. We will just say `left copure hereditary' when the class $\mathcal{Q}$ is known.

\begin{theorem} The following statements are equivalent:
\begin{enumerate}
\item The ring $R$ is left copure hereditary.
\item The homomorphic image of any copure Baer injective $R$-module is copure Baer injective.
\item The homomorphic image of any injective $R$-module is copure Baer injective.
\item Any finite sum of injective $R$-modules is copure Baer injective.
\end{enumerate}
\begin{proof}
$(1) \Rightarrow (2)$ Consider the diagram 
$$\begin{tikzcd}
0 \arrow{r} & I\arrow[hook]{r}\arrow{d}{f}
&R &\\
M \arrow{r}{g}&K \arrow{r} &0
\end{tikzcd}$$
of $R$-modules, where $I$ is a copure left ideal in $R$ and $M$ is a copure Baer injective module. Projectivity of $I$ gives the existence of a $\phi: I \rightarrow M$ such that $g\phi = f$. Copure Baer injectivity of $M$ gives a map $\phi^\prime :R \rightarrow M$ extending $\phi$, hence $g\phi^\prime$ extends $f$ and $K$ is copure Baer injective.
$(2) \Rightarrow (3)$ is trivial. $(3) \Rightarrow (1)$ Let $I$ be a copure left ideal of $R$ and consider the following diagram for a given $R$-module $M$:
$$\begin{tikzcd}
0 \arrow{r} & I\arrow[hook]{r}{\iota}\arrow{d}{f}
&R &\\
\E(M) \arrow{r}{g}&K \arrow{r} &0
\end{tikzcd}$$
where $\E(M)$ denotes the injective envelope of $M$. Since $K$ is copure Baer injective, there is a map $h:R \rightarrow K$ such that $\restr{h}{I} = f$. Projectivity of $R$ gives a $\sigma : R \rightarrow \E(M)$ such that $g\sigma = h$, i.e. $g\sigma \iota = h\iota =f$. This means $I$ is $\E(M)$-projective, i.e. $I$ is projective.
$(3) \Rightarrow (4)$ is clear. $(4) \Rightarrow (3)$ Similar to the proof of $(4) \Rightarrow (3)$ in \cite[Theorem 2.2]{Nada}.
\end{proof}
\end{theorem}

\section{\bf{Imbedding in Copure Baer Injective Modules}}
The main result of this section is the following:
\begin{theorem} \label{imbedding}
Let $\mathcal{Q}$ be a non-empty class of $R$-modules. Every module can be imbedded as a $\mathcal{Q}$-copure submodule in some $\mathcal{Q}$-copure Baer injective module.
\end{theorem}

We break the proof into three lemmas:

\begin{lemma} \label{imbeddinglemma} Every module can be imbedded in a copure Baer injective module.
\begin{proof} Given a module $A$, we want to show the existence of a copure Baer injective module that contains $A$ as a submodule. Consider the copure left ideals $I$ of $R$ and the set $\mathscr F$ of all maps $f:I \rightarrow A$. Thus, for any $f \in \mathscr F$ there is a pushout $B$ and a map $g:R \rightarrow B$ with $\restr{g}{I} = f$. The module $B$ may not be copure Baer injective, so put $A_0 = A$, $A_1 = B$ and repeat the above process with $A$ replaced by $A_1$ to give $A_2$ and $A_0 \subseteq A_1 \subseteq A_2$. Continuing in this manner, we get a sequence $A_0 \subseteq \cdots \subseteq A_n \subseteq A_{n+1} \subseteq \cdots$, for all $n \in \mathbb{N}$. Put $A_\omega = \bigcup A_n$. Now, for each nonlimit ordinal repeat the above process. If we get to a limit ordinal, say $l$, define $A_l = \bigcup \{A_s, s \, \textless \, l \}$. Let $t$ be the smallest ordinal with cardinality bigger than that of the ring $R$, i.e. $\abs{t} =\abs{R}^+$ (the successor cardinal of $\abs{R}$). For each $s \, \textless \, t$, we have $\abs{s} \thinspace \textless \thinspace \abs{t}$, $t$ is an initial ordinal and $A_t = \bigcup \{A_s, s \thinspace \textless \thinspace t \}$.  Now, $A_t$ is our copure Baer injective module. To see this, let $I$ be a copure left ideal of $R$ and $f: L \rightarrow A_t$ any map.  For each $r \in I$, let $s(r)$ be the smallest ordinal such that $f(r) \in A_{s(r)}$. Then $s(r) \thinspace \textless \thinspace t$ and $\abs{s(r)} \thinspace \textless \thinspace \abs{t} =\abs{R}^+$. Hence $\abs{s(r)} \leq \abs{R}$. Put $p = \sup \{s(r), r \in R\}$. As each $\abs{s(r)} \leq \abs{R}$, we must have $\abs{p} \leq \abs{R} \, \textless \,\, \abs{t}$. Hence, $p \thinspace \textless \thinspace t$. Since $t$ is a limit ordinal, we have $p+1 \thinspace \textless \thinspace t$. Therefore, for each $r \in R$, $r \in A_{s(r)} \subseteq A_p \subseteq A_{p+1} \subseteq A_t$. So, $f(I) \subseteq A_p$. Moreover, the map $f:I \rightarrow A_p$ can be extended to a map $g:R \rightarrow A_{p+1}$ with $\restr{g}{I} = f$. View $g$ now as a map $R \rightarrow A_t$. (The proof is adapted from \cite[p. 295]{Dauns}.)
\end{proof}
\end{lemma}

Of course, we know that every module can be imbedded in an injective (hence copure Baer injective) module. But this, unlike the next lemmas, does not guarantee that the imbedding is copure.

\begin{lemma} \label{chain}
Suppose that $A_0 \subseteq A_1 \subseteq \cdots$ is an ascending chain of modules such that $A_i$ is a copure submodule of $A_{i+1}$ for all $i$. Then, $A_0$ is copure in $\bigcup A_i$.
\begin{proof} Let $M$ be a member of the class $\mathcal{Q}$ and $f_0 : A_0 \rightarrow M$ a map which extends, by assumption, to a map $f_1 : A_1 \rightarrow M$, which in turn extends to $f_2 : A_2 \rightarrow M$, and so on. View the maps $f_i$ as sets of ordered pairs $(a_i,f(a_i))$ with $a_i \in A_i$ for all $i$. Hence, it is clear that $f_i \subseteq f_{i+1}$ for all $i$ and if $(x,y_1)$, $(x,y_2)$ $\in f_i$ for some $i$ then $y_1 = y_2$. Now, claim that $f = \cup f_i$ is a (well-defined) homomorphism. To see this, let $x \in \cup A_i$, i.e. $x \in A_i$ for some $i$ and $(x,f_i (x)) \in f_i \subseteq f$. If $(x,y_1), (x,y_2) \in f$, then $(x,y_1) \in f_i$ and $(x,y_2) \in f_j$ for some $i$ and $j$. Without loss of generality, assume $i \leq j$, so that $(x,y_1)$ and $(x,y_2)$ are both in $f_j$ and therefore $(x,y_1) = (x,y_2)$ and $f$ is well-defined. To finish the proof, let $x, y \in \cup A_i$ so that $x \in A_i$ and $y \in \cup A_j$ for some $i$ and $j$. Again assume $i\leq j$, so $f_j (x) = f_i(x)$. Now, for any $r, s \in R$, $f_j(rx+sy) = rf_j(x)+sf_j(y)$. So $f(rx+sy) = rf(x) + sf(y)$.
\end{proof}
\end{lemma}

\begin{lemma} The imbedding in Lemma \ref{imbeddinglemma} is copure.
\begin{proof} The construction of $A_i$ in Lemma \ref{imbeddinglemma} shows, by (1) of  Lemma \ref{pushout}, that $A_i$ is copure in $A_{i+1}$ for all $i$, and by Lemma \ref{chain}, $A$ is copure in $\bigcup A_n = A_\omega$. Again by Lemma \ref{chain}, $A_\omega$ is copure in $A_{\omega + 1}$ and $A_{\omega + 1}$ is copure in $A_{\omega + 2}$ and so on. In other words, for every ordinal $s \thinspace \textless \thinspace \abs{R}^+$, we have either $A$ is copure in $A_s$ if $s$ is not a limit ordinal, or $A_s = \bigcup_{u \textless s} A_u$ if $s$ is a limit ordinal. In either case, $A$ is copure in $A_t$, as desired.
\end{proof}
\end{lemma}

The imbedding Theorem can be used in characterizing some copure exact sequences.

\begin{theorem} \label{cpbiinj1}The sequence $0 \rightarrow I \stackrel{\iota}{\hookrightarrow} R \rightarrow R/I \rightarrow 0$ is copure exact if and only if every copure Baer injective $R$-module is injective with respect to it.
\begin{proof} Necessity is clear. To prove sufficiency, let $j:I \rightarrow C$ be a copure imbedding in a copure Baer injective module $C$ (Theorem \ref{imbedding}). Therefore, there exists a map $f:R \rightarrow C$ such that $f \iota = j$. But $j$ is a copure monomorphism, so by (2) of Lemma \ref{pushout} $\iota$ is a copure monomorphism.
\end{proof}
\end{theorem}

\section{\bf{Characterization of Rings Using Copure Baer Injectivity}}
Thani \cite{Nada} proved that for a left self injective ring $R$, the condition that $R/I$ is pure Baer injective for every essential left ideal $I$ of $R$ is enough to make $R/I$ pure Baer injective for all left ideals $I$ of $R$. Using the same line of argument, we  generalize this to copure Baer injectivity.
\begin{proposition}
Let $R$ be a left self injective ring. If $R/J$ is copure Baer injective for any essential left ideal $J$ of $R$, then $R/I$ is copure Baer injective for any left ideal $I$ of $R$.
\begin{proof} Since $R$ is injective, the injective envelope $\E(I)$ of $I$ must be a direct summand of $R$, for any left ideal $I$ of $R$. Therefore, $\E(I) = Re$ for some idempotent $e \in \E(I)$. Now for the map $f:R \rightarrow Re$, defined by $1 \mapsto e$, since $I$ is essential in $R$, $f^{-1}(I)$ must as well be essential in $R$ and, therefore by assumption, $R/f^{-1}(I)$ is copure Baer injective. Define $\xbar{f}: Re/I \rightarrow R/f^{-1}(I)$ by $re+I \mapsto r+f^{-1}(I)$ and proceed as in the proof of \cite[Proposition 2.3]{Nada}.
\end{proof}
\end{proposition}

By a $\mathcal{Q}$-\emph{copure split} ring we mean a ring every $\mathcal{Q}$-copure left ideal of which is a direct summand (=principal ideal). Clearly, every pure split ring is $\mathcal{PI}$-copure split and if a ring $R$ is left copure-split then it is left copure hereditary. The $\mathcal{Q}$-copure split rings are characterized in the following Theorem.

\begin{theorem} \label{cpsplit} The following statements are equivalent:
\begin{enumerate}
\item The ring $R$ is left copure split.
\item Every $R$-module is copure Baer injective.
\item Any copure left ideal of $R$ is copure Baer injective.
\item \begin{enumerate}
\item $R$ is left copure hereditary, and
\item Every free left $R$-module is copure Baer injective.
\end{enumerate}
\end{enumerate}
\begin{proof} $(1) \Rightarrow (2)$ Let $M$ be an $R$-module. Since every left ideal $I$ of $R$ is a direct summand, every map $I \rightarrow M$ into any $R$-module can easily be extended  to a map $R \rightarrow M$. $(2) \Rightarrow (3)$ is obvious. $(3) \Rightarrow (1)$ Let $I$ be a copure left ideal of $R$. Copure Baer injectivity of $I$ gives a homomorphism $R \rightarrow I$ that extends the identity map of $I$, which means $I$ is a direct summand of $R$. $(1) \Rightarrow (4)$(a) and $(2) \Rightarrow (4)$(b) are immediate. $(4) \Rightarrow (3)$ Let $I$ be a copure left ideal of $R$, hence projective by (a) and therefore a direct summand of some free $R$-module $F$. From (b) it follows that $F$ is copure Baer injective and by Proposition \ref{Ext}, so is $I$.
\end{proof}
\end{theorem}

Recall that a ring $R$ is called \emph{left coregular} if every left ideal of $R$ is copure in $R$ \cite[p.324]{Wis}.

\begin{theorem} \label{cpbiinj} For a ring $R$ the following statements are equivalent:
\begin{enumerate}
\item The ring $R$ is left coregular.
\item Every copure Baer injective $R$-module is injective.
\item Every copure Baer injective $R$-module is quasi injective.
\end{enumerate}
\begin{proof} $(1) \Rightarrow (2) \Rightarrow (3)$ are obvious. $(2) \Rightarrow (1)$ By assumption, every copure Baer injective $R$-module is injective with respect to any sequence $0 \rightarrow I \rightarrow R \rightarrow R/I \rightarrow 0$, which must, therefore, be copure exact by Theorem \ref{cpbiinj1}. $(3) \Rightarrow (2)$  Let $M$ be a copure Baer injective $R$-module. Hence, by Proposition \ref{Ext}, so is $M \oplus \E(R)$ which must be quasi injective by assumption. Therefore, $M$ is injective with respect to $\E(R)$. In particular, $M$ is $R$-injective or injective by Baer condition.
\end{proof}
\end{theorem}

\bibliographystyle{amsplain}

\end{document}